\documentclass{ws-rv961x669}
\usepackage{ws-rv-van}     
\usepackage{ws-rv-thm}     
\usepackage{subfigure}     
\makeindex
\usepackage{amsmath}
\usepackage{amssymb}
\usepackage{ulem,bm}

\begin{document}

\chapter[Model Reduction in Stochastic Environments]{Model Reduction in Stochastic Environments}\label{ra_ch1}

\author[E. Forgoston, L. Billings \& I.B. Schwartz]{Eric
  Forgoston,\footnote{eric.forgoston@montclair.edu} and Lora Billings}

\address{Department of Mathematical Sciences,\\ Montclair State University,
  Montclair, NJ 07043, USA. \\}

\author[E. Forgoston, L. Billings \& I.B. Schwartz]{Ira B. Schwartz}

\address{Nonlinear Systems Dynamics Section, Plasma Physics Division, Code 6792,\\ U.S. Naval Research Laboratory, Washington, DC 20375, USA
 \\}

\begin{abstract}
We present a general
theory of stochastic model reduction which is based on a normal form
coordinate transform method of A.J. Roberts. This nonlinear, stochastic projection allows for the
deterministic and stochastic dynamics to interact correctly on the
lower-dimensional manifold so that the dynamics predicted by the reduced,
stochastic system agrees well with the dynamics predicted by the original,
high-dimensional stochastic system. The method may be applied to any
system with well-separated time scales. In this article, we
consider a physical problem that involves a singularly perturbed Duffing oscillator as well as a biological problem that involves the
prediction of infectious disease outbreaks.
\end{abstract}


\body


\section{Introduction}\label{sec:intro}
It is well-known that noise can have a significant effect on deterministic
dynamical systems. As an example, given an initial state starting in a basin
of attraction, noise can cause the initial state to cross the basin boundary
and move into another, distinct basin of
attraction~\cite{dyk90,dmsss92,mil96,lumcdy98,forsch09}.  Many researchers have 
investigated how noise affects physical and biological phenomena at a wide
variety of levels including
switching between the magnetisation states in magnets~\cite{kohn2005}, and voltage and current
states in Josephson junctions~\cite{Fulton1974}, sub-cellular processes, tissue dynamics, large-scale population
dynamics~\cite{Tsimring_2014}, genetic switching~\cite{Assaf2011}, and extinction and switching in general heterogeneous networks~\cite{lindley2014rare,schwartz2016epidemic,hindes2017epidemic,hindes2017large}. 

Stochasticity manifests itself as either external or internal noise. In this
article, we shall consider only external noise, which comes from a source outside the system being considered (e.g. population
growth under the influence of climatic effects, or a random signal fed into a
transmission line), and often is modeled by replacing an
external parameter with a random process. Mathematically, the effect of external noise is often described using a Langevin
equation or the associated Fokker-Planck equation (though the dynamics of
external noise may sometimes be described by a master equation~\cite{Roberts2015}).

The reduction in high-dimensional systems is an important and fundamental
problem in nonlinear dynamical systems. Moreover, normal form coordinate
transforms provide a way to simplify multiscale nonlinear dynamics by
separating the long-term dynamics of interest from the transient dynamics~\cite{rob08}.
In this article, we present a general
theory of stochastic model reduction that is based on a normal form
coordinate transform method of A.J. Roberts. This nonlinear, stochastic projection allows for the
deterministic and stochastic dynamics to interact correctly on the
lower-dimensional manifold so that the transformed dynamics reproduce
fully the original dynamics. The method may be applied to any
system with well-separated time scales. Here, we
consider a physical problem that is associated with the use of unmanned sensors
operating in the ocean as well as a biological problem that involves the
prediction of infectious disease outbreaks.

Section~\ref{sec:gen_theory} provides an overview of the theory including the
use of center manifold theory for deterministic problems (Sec.~\ref{sec:gen_theory_DCM}) and the normal form
coordinate transform for stochastic problems
(Sec.~\ref{sec:gen_theory_SCM}). Sections~\ref{sec:ex1} and~\ref{sec:ex2}
provide two examples - the first involves a singularly perturbed
stochastic Duffing oscillator, while the second involves a stochastic
Susceptible-Exposed-Infectious-Recovered (SEIR) epidemic model.  In the first
example, we show that one can use deterministic theory as the noise effects
occur at such high order so that the stochastic correction is negligible. The
reduced system is used to understand a variety of system behaviour including
the optimal escape path and escape rate. In the second example, we show how
one must use the stochastic theory to obtain long-time predictions of disease
outbreak. Conclusions are found in Sec.~\ref{sec:conc}.

\section{General Theory}\label{sec:gen_theory}
We consider the following general $(m+n)$-dimensional system of Stratonovich stochastic differential equations
\begin{subequations}
\begin{equation}
\label{e:Gen_sys_x}
\dot{\bf{x}}={\bf A}{\bf x} + {\bf F}\left ({\bf x},{\bf y},{\bm \Phi}\right ),
\end{equation}  
\begin{equation}
\label{e:Gen_sys_y}
\dot{\bf{y}}={\bf B}{\bf y} + {\bf G}\left ({\bf x},{\bf y},{\bm \Psi}\right ),
\end{equation}  
\end{subequations}
where  ${\bf x}(t)\in\mathbb{R}^m$, ${\bf
  y}(t)\in\mathbb{R}^n$, ${\bm \Phi}(t)$ and ${\bm \Psi}(t)$ describe stochastic forces with
adjustable noise intensity, ${\bf A}$ and ${\bf B}$ are constant matrices, and
${\bf F}$ and ${\bf G}$ are stochastic, nonlinear functions.
When there exist slow and fast time scales, the special case of stochastic
singular perturbation systems have been explored for
  realisations~\cite{bergen03,berglund2006noise} and for probabilistic large
fluctuations~\cite{heckman2014stochastic}, and have the form
\begin{subequations}
\begin{equation}
\label{e:Gen_sys_x-p}
\dot{\bf{x}}={\bf A}{\bf x} + {\bf F}\left ({\bf x},{\bf y},{\bm \Phi}\right ),
\end{equation}  
\begin{equation}
\label{e:Gen_sys_y-p}
\epsilon\,\dot{\bf{y}}={\bf B}{\bf y} + {\bf G}\left ({\bf x},{\bf y},{\bm \Psi}\right ),
\end{equation}  
\end{subequations}
where $\epsilon$ is a small parameter.

\subsection{Deterministic Center Manifold}\label{sec:gen_theory_DCM}

To begin, we remove the stochastic terms from Eqs.~(\ref{e:Gen_sys_x})-(\ref{e:Gen_sys_y}) so that ${\bf F}={\bf F}({\bf x},{\bf y})$ and ${\bf
  G}={\bf G}({\bf x},{\bf y})$.  
A general nonlinear system may be transformed so that the system's
linear part has a block diagonal form consisting of three matrix blocks.  The
first matrix block will possess eigenvalues with positive real part; the
second matrix block will possess eigenvalues with negative real part; and the
third matrix block will possess eigenvalues with zero real part.  These three
matrix blocks are respectively associated with the unstable eigenspace, the
stable eigenspace, and the center eigenspace.  If there are no eigenvalues
with positive real part, then the orbits will rapidly decay to the center
eigenspace.

It is often the case that a system of equations can not be written in a block
diagonal form with one matrix block possessing eigenvalues with negative real
part and the other matrix block possessing eigenvalues with zero real
part. Even though it is possible to construct a center manifold from a system not in separated block form~\cite{chilat97}, it is much easier to apply the center
  manifold theory to a system with separated stable and center directions.
Therefore, we generally transform the original system of equations to a new system of equations that will have
the eigenvalue structure that is needed to apply rigorous center manifold theory~\cite{car81}.  The
theory allows one to find an invariant center manifold that passes through a
fixed point and to which one can restrict the new transformed system. 

To make the ideas more concrete, consider the singularly perturbed system
given by Eqs.~(\ref{e:Gen_sys_x-p})-(\ref{e:Gen_sys_y-p}), and let $t=\epsilon\tau$.  Denoting $\dot{}$ as
$d/dt$ and $^{\prime}$ as $d/d\tau$, then the deterministic form of
Eqs.~(\ref{e:Gen_sys_x-p})-(\ref{e:Gen_sys_y-p}) is transformed to the
following system of equations: 
\begin{subequations}
\begin{equation}
\label{e:trans_sys_x}
{\bf{x}}^{\prime}=\epsilon\left ({\bf A}{\bf x} + {\bf F}({\bf x},{\bf
    y})\right ),
\end{equation}  
\begin{equation}
\label{e:trans_sys_y}
{\bf{y}}^{\prime}={\bf B}{\bf y} + {\bf G}({\bf x},{\bf y}),
\end{equation}
\begin{equation}
\label{e:trans_sys_eps}
\epsilon^{\prime}=0.
\end{equation}  
\end{subequations}

We recast the problem by treating $\epsilon$ as
  a state variable, and we let ${\bf \bar{A}}=\epsilon {\bf A}$ and ${\bf
    \bar{F}}=\epsilon {\bf
    F}$. Equations~(\ref{e:trans_sys_x})-(\ref{e:trans_sys_eps}) can thus be
  rewritten as
\begin{subequations}
\begin{equation}
\label{e:rewrit_x}
{\bf{x}}^{\prime}={\bf \bar{A}}{\bf x} + {\bf \bar{F}}({\bf x},{\bf y},\epsilon),
\end{equation}
\begin{equation}
\label{e:rewrit_y}
{\bf{y}}^{\prime}={\bf B}{\bf y} + {\bf G}({\bf x},{\bf y}),
\end{equation}
\begin{equation}
\label{e:rewrit_eps}
\epsilon^{\prime}=0.
\end{equation}  
\end{subequations}
If ${\bf \bar{A}}$ and ${\bf B}$ are constant matrices such that all
of the eigenvalues of ${\bf \bar{A}}$ have zero real parts, while all of the
eigenvalues of ${\bf B}$ have negative real parts, then the system will
rapidly collapse onto a lower-dimensional manifold given by center manifold
theory~\cite{car81}.  

If the center manifold is assumed to be smooth and given by 
\begin{equation}
\label{e:cent_man}
{\bf y}={\bf h}({\bf x},\epsilon),
\end{equation}
then substitution of Eq.~(\ref{e:cent_man}) into Eq.~(\ref{e:rewrit_y}) leads
to the following center manifold condition:
\begin{equation}
\label{e:cent_man_cond}
{\bf h}_{{\bf x}}\left ({\bf \bar{A}}{\bf x} + {\bf \bar{F}}({\bf
  x},{\bf h}({\bf x},\epsilon),\epsilon)\right ) = {\bf B}{\bf h}({\bf
  x},\epsilon) + {\bf G}({\bf
  x},{\bf h}({\bf x},\epsilon)),
\end{equation}
where ${\bf h}_{{\bf x}}$ denotes the partial derivative  of ${\bf h}$ with
respect to ${\bf x}$.  Although it is
generally not possible to solve Eq.~(\ref{e:cent_man_cond}) for ${\bf h}$, one can
approximate the center manifold by expanding ${\bf h}$ in the following way:
\begin{equation}
\label{e:expand}
{\bf h}({\bf x},\epsilon)={\bf h_0}({\bf x})+\epsilon {\bf
  h_1}({\bf x})+\epsilon^2 {\bf h_2}({\bf x})+\mathcal{O}(\epsilon^3).
\end{equation}
Typically, this approximation of ${\bf h}({\bf x},\epsilon)$ is found by
substituting Eq.~(\ref{e:expand}) into the center manifold condition
(Eq.~(\ref{e:cent_man_cond})) and matching coefficients.

\subsection{Stochastic Center Manifold and the Normal Form Coordinate Transform}\label{sec:gen_theory_SCM}

In a manner similar to that shown in Sec.~\ref{sec:gen_theory_DCM}, the
singularly perturbed stochastic system given by Eqs.~(\ref{e:Gen_sys_x-p})-(\ref{e:Gen_sys_y-p}) can be transformed~\cite{bergen03} to the form
given by Eqs.~(\ref{e:rewrit_x})-(\ref{e:rewrit_eps}), where now ${\bf
  \bar{F}}={\bf \bar{F}}({\bf x},{\bf y},\epsilon,{\bm \Phi})$ and ${\bf
  G}={\bf G}({\bf x},{\bf y},{\bm \Psi})$.  If ${\bf
  \bar{A}}$ and ${\bf B}$ satisfy the same spectral conditions as for
the deterministic system, and if the stochastic time dependence found in ${\bf
  \bar{F}}$ and ${\bf
  G}$ is due to independent white noise processes, then there exists a
stochastic center manifold for the original stochastic system~\cite{box89}. 

One method for computing the stochastic center manifold for systems with
both fast and slow dynamics uses the construction of a normal form coordinate
transform that not only reduces the dimension of the dynamics, but also
separates all of the fast processes from all of the slow
processes~\cite{rob08}.  While this type of normal form coordinate transform
may be used to find deterministic center manifolds, the application of this
transform to stochastic systems is particularly interesting since white noise
has fluctuations on all scales that can lead to unbounded solutions.

There are many publications~\cite{knowie83,coelti85,nam90,namlin91} which deal with the simplification of a stochastic
dynamical system using a stochastic normal form transformation.  In these
articles, the noise term is multiplied by a small parameter, and therefore,
the resulting stochastic normal form is a perturbation of the deterministic
normal form.  Furthermore, one can find in Coullet et al.,~\cite{coelti85} and
Namachchivaya et al.,~\cite{namlin91} normal
form transformations that involve anticipative noise processes.  However,
these integrals of the noise process into the future were not dealt with
rigorously.

Rigorous, theoretical analysis to support normal form coordinate transforms
(and center manifold reduction) was developed by Arnold and
Imkeller~\cite{arnimk98} and Arnold~\cite{arn98}, where the technical problem of the anticipative noise integrals also was
dealt with rigorously.  Later, another stochastic normal form transformation
was developed by Roberts~\cite{rob08}.  This new method is such that ``anticipation can ... always [be] removed from
the slow modes with the result that no anticipation is required after the fast
transients decay''(Roberts~\cite{rob08}, pp. 13).  An advantage of removing
anticipation is the simplification of the normal form.  Nonetheless, this
simpler normal form retains its accuracy with the original stochastic system.
Furthermore, when modeling the macroscopic behaviour of microscopic,
stochastic systems, it is desirable to avoid anticipation in the normal form~\cite{rob08}. 
It is important to note that the normal form is valid for all time since it is
just a coordinate transform.  Furthermore, the dynamics also are valid for all
time as long as the truncation error is small enough for the problem of
interest.

In the examples presented in Sec.~\ref{sec:ex1} and Sec.~\ref{sec:ex2}, we
shall use the method of Roberts~\cite{rob08} to
simplify our stochastic dynamical system to one that emulates the long-term dynamics of
the original, multiple-time-scale system.  The method involves five
principles, which we recapitulate here for the purpose of clarity.  The
principles are as follows:
\begin{enumerate}
\item Avoid unbounded, secular terms in both the transformation and the
  evolution equations to ensure a uniform asymptotic approximation.
\item Decouple all of the slow processes from the fast processes to ensure a
  valid long-term model.
\item Insist that the stochastic slow manifold is precisely the transformed
  fast processes coordinate being equal to zero.
\item To simplify matters, eliminate as many as possible of the terms in the evolution equations.
\item Try to remove all fast processes from the slow processes by avoiding as
  much as possible the fast time memory integrals in the evolution equations.
\end{enumerate}

In practice, the original stochastic system of equations (which satisfy the
necessary spectral requirements) in $({\bf x},{\bf y})$ coordinates is
transformed to a new $({\bf X},{\bf Y})$ coordinate system using a stochastic
coordinate transform as follows:
\begin{subequations}
\begin{equation}
\label{e:stoch_coord_trans_x}
{\bf x}={\bf X}+{\bm \xi}({\bf X},{\bf Y},t),
\end{equation}
\begin{equation}
\label{e:stoch_coord_trans_y}
{\bf y}={\bf Y}+{\bm \eta}({\bf X},{\bf Y},t),
\end{equation}
\end{subequations}
where the specific form of Eqs.~(\ref{e:stoch_coord_trans_x}) and~(\ref{e:stoch_coord_trans_y}) is chosen to simplify the original
system according to the five principles listed previously.  The terms ${\bm
  \xi}({\bf X},{\bf Y},t)$ and ${\bm \eta}({\bf X},{\bf Y},t)$ are found using
an iterative procedure that will be explicitly demonstrated using the first
example (Sec.~\ref{sec:ex1}) of a singularly
perturbed, damped, stochastic Duffing oscillator model.  Theoretical
details of the normal form coordinate transform process can be found in Roberts~\cite{rob08}.

\section{Example 1: Singularly Perturbed Stochastic Duffing Oscillator}\label{sec:ex1}
An important application in many fields is that of sensing
in stochastic environments.  Improved environmental sensing and prediction can
be achieved 
through the incorporation of continuous
monitoring of the region of interest.  For example, one could monitor the
stochastic ocean using autonomous underwater
gliders~\cite{webbsijo01,shermandov01,eriksenolwlsbc01}.  However, to do this,
one must understand both the dynamics and control of the gliders.

Extending the lifetime (energy optimisation problem) of sensing devices
(e.g. gliders) in stochastic environments such as the ocean requires
an understanding of the effect of the environmental forces on both the devices
and the
region being monitored.  The ocean dynamics are high-dimensional and stochastic.  Therefore, as a first step towards using the underlying ocean structure to optimize a
sensor's energy usage, 
we will use the two
methods described in Sec.~\ref{sec:gen_theory} to obtain a reduction in the dimension of
the stochastic system.  The reduced system can then be used to understand a
variety of system behaviour including the optimal escape path and the escape rate.

We consider the following singularly perturbed, damped, Duffing oscillator
system with additive noise (see Heckman and Schwartz~\cite{heckman2014stochastic} for a large
fluctuation approach that is global):
\begin{subequations}
\begin{equation}
\label{e:Duff_sys_x}
\dot{x}=y + \sqrt{2D}\phi(t),
\end{equation}
\begin{equation}
\label{e:Duff_sys_y}
\epsilon\dot{y}=(x-x^3-y),
\end{equation}
\end{subequations}
where $D$ is the noise intensity and $\phi(t)$ describes a stochastic white
force that is characterized by the following correlation
functions:
\begin{subequations}
\begin{equation}
\label{e:mean}
\langle\phi_i(t)\rangle=0,
\end{equation}
\begin{equation}
\label{e:corr}
\langle\phi_i(t)\phi_j(t^{\prime})\rangle=\delta(t-t^{\prime})\delta_{ij}.
\end{equation}
\end{subequations}

\subsection{Deterministic Center Manifold}\label{sec:example_DCM}

Following the general theory of Sec.~\ref{sec:gen_theory_DCM}, we consider the deterministic form
of Eqs.~(\ref{e:Duff_sys_x})-(\ref{e:Duff_sys_y}) by setting $\phi(t)=0$.  The slow manifold is found by setting $\epsilon =0$ in Eq.~(\ref{e:Duff_sys_y}).
Solving for $y$ gives the equation of the slow manifold as $y=x-x^3$ (which
corresponds to $h_0(x)$ in Eq.~(\ref{e:expand})).
Substitution of this into the deterministic form of Eq.~(\ref{e:Duff_sys_x}) gives the dynamics along the slow
manifold as $\dot{x}=x-x^3$.

If, as in Sec.~\ref{sec:gen_theory_DCM}, we let $t=\epsilon\tau$ and denote $\dot{}$ as
$d/dt$ and $^{\prime}$ as $d/d\tau$, then Eqs.~(\ref{e:Duff_sys_x})-(\ref{e:Duff_sys_y}) (with $\phi(t)=0$)
are transformed to the following system:
\begin{subequations}
\begin{equation}
\label{e:Duff_trans_x}
x^{\prime}=\epsilon y, 
\end{equation}
\begin{equation}
\label{e:Duff_trans_y}
y^{\prime}=x-x^3-y, 
\end{equation}
\begin{equation}
\label{e:Duff_trans_eps}
\epsilon^{\prime}=0.
\end{equation}
\end{subequations}
Rearrangement of Eqs.~(\ref{e:Duff_trans_x})-(\ref{e:Duff_trans_eps}) leads to
a system described by constant matrices ${\bf \bar{A}}$ and ${\bf B}$ that satisfy the spectral
requirements of Sec.~\ref{sec:gen_theory_DCM}.  Furthermore, since the $x$ and $\epsilon$
variables are associated with the ${\bf \bar{A}}$ matrix (eigenvalues with
zero real parts), and
the $y$ variable is associated with the ${\bf B}$ matrix (eigenvalues with
negative real parts), we know that the
center manifold is given by $y=h(x,\epsilon)$.

The center manifold condition is given by Eq.~(\ref{e:cent_man_cond}), and we
approximate the center manifold (Eq.~(\ref{e:expand})) as follows:
\begin{subequations}
\begin{flalign}
h(x,\epsilon)&=h_0(x)+\epsilon h_1(x)+\epsilon^2 h_2(x)+\mathcal{O}(\epsilon^3)\label{e:h}\\
&=c_0 + c_{01}\epsilon +c_{10}x +c_{02}\epsilon^2 + c_{11}x\epsilon +
c_{20}x^2 \nonumber\\
&\hspace{0.77cm}+ c_{03}\epsilon^3 + c_{12}x\epsilon^2 + c_{21}x^2\epsilon +
c_{30}x^3 + \mathcal{O}(\gamma^4),\label{e:h2}
\end{flalign}
\end{subequations}
where $c_0$, $c_{01}$, $c_{10}$, $c_{02}$, $\ldots$ are unknown coefficients,
and $\gamma=|(x,\epsilon)|$ so that $\gamma$ provides a count of the number of
$x$ and $\epsilon$ factors in any one term.
The center manifold condition for this example
is given by
\begin{equation}
\label{e:Duff_cent_man_cond}
\frac{\partial h(x,\epsilon)}{\partial x}(\epsilon h(x,\epsilon))=-h(x,\epsilon)+x-x^3.
\end{equation}

By substituting Eq.~(\ref{e:h2}) into Eq.~(\ref{e:Duff_cent_man_cond}) and matching the
different orders to find the coefficients, one finds the following center
manifold equation (expanded to sixth-order):
\begin{subequations}
\begin{flalign}
h(x,\epsilon) = 
&\, \, x-x^3+\epsilon (-x+4x^3-3x^5) + \epsilon^2 (2x-20x^3)\nonumber\\
&\hspace{0.31cm} +\epsilon^3
(-5x+104x^3)+\epsilon^4 (14x) +\epsilon^5 (-42x)+\mathcal{O}(\gamma^7).\label{e:h4}
\end{flalign}
\end{subequations}
Note that by letting $\epsilon =0$, one recovers the zero-order approximation,
$h_0(x)$ (the slow manifold).  In addition, since $\epsilon$ is now a state
variable, the first nontrivial correction term to the zero-order approximation
is a quadratic term.

 \subsection{Stochastic Center Manifold and the Normal Form Coordinate Transform}\label{sec:example_SCM}

To describe the stochastic effects, we will derive the
normal form coordinate transform (and thus the stochastic center manifold) for the singularly perturbed, stochastic Duffing system given by
Eqs.~(\ref{e:Duff_sys_x})-(\ref{e:Duff_sys_y}).  As demonstrated previously, use of
the $t=\epsilon \tau$ transformation leads to the following system:
\begin{subequations}
\begin{equation}
\label{e:Duff_trans_x_stoch}
x^{\prime}=\epsilon (y +\sqrt{2D}\phi)=\epsilon(y+\sigma\phi),
\end{equation}
\begin{equation}
\label{e:Duff_trans_y_stoch}
y^{\prime}=x-x^3-y,
\end{equation}
\begin{equation}
\label{e:Duff_trans_eps_stoch}
\epsilon^{\prime}=0,
\end{equation}
\end{subequations}  
where $\sigma$ is the
standard deviation of the noise intensity $D=\sigma^2/2$.

The construction of the
normal form is quite tedious and complicated (although it is possible to
derive the normal form using a computer algebra system).  However, the result allows one
to determine if there are any noise terms that cause a significant difference between the
average stochastic center manifold (the stochastic center manifold generally
fluctuates about an average location) and the deterministic center manifold.

For this problem, it turns out that the noise terms that could lead to a
difference between the deterministic and average stochastic center manifolds occur at very high order in the
normal form expansion.  Therefore, the correction to the deterministic center
manifold is minimal, and we expect that one can use the deterministic center
manifold results of
Sec.~\ref{sec:example_DCM} to accurately solve the problem of interest. In
our case, we are interested in understanding the optimal escape path and the
escape rate, and we expect that results found using the deterministic center
manifold reduction will agree very well with numerical computations using the
original stochastic system (Eqs.~(\ref{e:Duff_sys_x})-(\ref{e:Duff_sys_y})).

We proceed by showing how to use the method
of Roberts~\cite{rob08} described in Sec.~\ref{sec:gen_theory_SCM} to construct a normal form coordinate
transform that separates the slow and fast dynamics of
Eqs.~(\ref{e:Duff_trans_x_stoch})-(\ref{e:Duff_trans_y_stoch}).  In what
follows, we outline the steps involved in the first iteration, while details
regarding the higher iterations can be found in Forgoston and Schwartz~\cite{forsch09}.

\subsubsection{First Iteration}\label{example_SCM_first}

We begin by letting
\begin{subequations}
\begin{equation}
\label{e:xNF}
x\approx X,
\end{equation}
\begin{equation}
\label{e:XdotNF}
X^{\prime}\approx 0,
\end{equation}
\end{subequations}
and by finding a change to the $y$ coordinate (fast process) with the form
\begin{subequations}
\begin{equation}
\label{e:yNF}
y = Y+\eta(\tau,X,Y)+\ldots,
\end{equation}
\begin{equation}
\label{e:YdotNF}
Y^{\prime} = -Y+G(\tau,X,Y)+\ldots,
\end{equation}
\end{subequations}
where $\eta$ and $G$ are small corrections to the coordinate transform and the
corresponding evolution equation.  Substitution of
Eqs.~(\ref{e:xNF})-(\ref{e:YdotNF}) into Eq.~(\ref{e:Duff_trans_y_stoch})
gives the following equation:
\begin{equation}
Y^{\prime}+\frac{\partial \eta}{\partial \tau} + \frac{\partial \eta}{\partial
  X}\frac{\partial X}{\partial \tau} + \frac{\partial \eta}{\partial
  Y}\frac{\partial Y}{\partial \tau} = -Y -\eta +X-X^3.
\end{equation}
Replacing $Y^{\prime}=\partial Y/\partial \tau$ with $-Y+G$ (Eq.~(\ref{e:YdotNF})),
noting that $\partial X/\partial \tau =0$ (Eq. (\ref{e:XdotNF})), and ignoring the
term $\partial\eta/\partial Y\cdot G$ since it is a product of small
corrections leads to the following:
\begin{equation}
\label{e:Getaeqn}
G+\frac{\partial \eta}{\partial \tau} -Y \frac{\partial \eta}{\partial
  Y}+\eta = X-X^3.
\end{equation}

Equation~(\ref{e:Getaeqn}) must now be solved for $G$ and $\eta$.  In order to
keep the evolution equation (Eq.~(\ref{e:YdotNF})) as simple as possible
(principle~(4) of Sec.~\ref{sec:gen_theory_SCM}), we let
$G=0$, which means that the coordinate transform (Eq.~(\ref{e:yNF})) is modified by $\eta =X-X^3$.
Therefore, the new approximation of the coordinate transform and its dynamics
are given by
\begin{subequations}
\begin{equation}
\label{e:yNF2_a}
y = Y+X-X^3+\mathcal{O}(\zeta^2),
\end{equation}
\begin{equation}
\label{e:yNF2_b}
Y^{\prime} = -Y+\mathcal{O}(\zeta^2),
\end{equation}
\end{subequations}
where $\zeta =|(X,Y,\epsilon,\sigma)|$ so that $\zeta$ provides a count of the
number of $X$, $Y$, $\epsilon$, and $\sigma$ factors in any one term. 

\subsubsection{Higher Iterations}\label{sec:example_SCM_higher}

The construction of the normal form continues by seeking corrections, $\xi$ and $F$, to the
$x$ coordinate transform and the $X$ evolution using the
updated residual of the $x$ equation (Eq.~(\ref{e:Duff_trans_x_stoch})), and
by seeking corrections, $\eta$ and $G$, to the
$y$ coordinate transform and the $Y$ evolution equation using the
updated residual of the $y$ equation (Eq.~(\ref{e:Duff_trans_y_stoch})).
Details regarding the second iteration can be found in Forgoston and Schwartz~\cite{forsch09}.

The
derivation of $\xi$ and $F$ in the second and fourth iterations along with the derivation
of $\eta$ and $G$ in the third iteration leads to the following updated approximation of the coordinate transforms and
their corresponding evolution equations:
\begin{subequations}
\begin{flalign}
y=&Y+X-X^3+\epsilon\left (-X+4X^3-3X^5\right )\nonumber\\
&\hspace{0.28cm}+\epsilon\sigma\left (-e^{-\tau}*\phi+3X^2e^{-\tau}*\phi\right )+3\epsilon^2XY^2+\mathcal{O}(\zeta^3),\label{e:yNF3_a}\\
Y^{\prime}=&-Y+\epsilon\left (-Y+3X^2Y\right )+\mathcal{O}(\zeta^3),\label{e:yNF3_b}\\
x=&X-\epsilon Y+\epsilon^2\left (Y-3X^2Y\right )\nonumber\\
&\hspace{0.32cm}+\epsilon^2\sigma\left (e^{-\tau}*\phi-3X^2e^{-\tau}*\phi\right )+\mathcal{O}(\zeta^4),\label{e:xNF4_a}\\
X^{\prime}=&\epsilon\left (X-X^3\right )+\epsilon\sigma\phi+\epsilon^2\left
  (-X+4X^3-3X^5\right )\nonumber\\
&\hspace{0.3cm}+\epsilon^2\sigma\left (-\phi+3X^2\phi\right )+\mathcal{O}(\zeta^4),\label{e:xNF4_b}
\end{flalign}
\end{subequations}
where the noise convolution is defined by
\begin{equation}
e^{-\tau}*\phi = \int\limits_{-\infty}^\tau \exp{\left [-(\tau-s)\right ]}\phi(s) \, ds.
\end{equation}
Details regarding the derivation of Eqs.~(\ref{e:yNF3_a})-(\ref{e:xNF4_b})
can be found in Forgoston and Schwartz~\cite{forsch09}.

One can continue this iterative procedure to obtain higher order terms in the
expansions of the coordinate transform and normal form.  For the stochastic
Duffing system under consideration, the fifth and sixth iterations lead to
updated approximations of the $x$ and $y$ coordinate transforms (along with
their associated evolution equations) that are extremely long and
complicated.  These approximations can be found in Forgoston and Schwartz~\cite{forsch09}.

In the higher order transform found after six iterations,
one can see the appearance of quadratic noise terms.  For example, one can see
terms of the form $e^{-\tau}*(e^{-\tau}*\phi)^2$ in the coordinate transforms, and one can see terms of
the form $\phi e^{-\tau}*\phi$ in one of the evolution equations.  This quadratic noise is
important because it leads to the creation of a deterministic drift within the
slow dynamics~\cite{rob08,namlin91}.  Furthermore, the stochastic center
manifold generally undergoes fluctuations about a mean or average location.
This average stochastic center manifold is usually different from the
deterministic center manifold, and it is the quadratic noise process that
generates this difference.

\subsubsection{Comparison with Deterministic Center Manifold and Effect of
  Quadratic Noise}\label{sec:example_SCM_comparison}

Letting $Y=0$ and $\sigma =0$ in the higher order transform found after six
iterations leads to the following deterministic center manifold equation:
\begin{subequations}
\begin{flalign}
\label{e:dCM_x}
x=&X,\\
\nonumber\\
y=&X-X^3+\epsilon(-X+4x^3-3X^5)+\epsilon^2(2X-20X^3+42X^5)\nonumber\\
&+\epsilon^3(-X+16X^3-66X^5+96X^7-45X^9)+\mathcal{O}(\epsilon^3).\label{e:dCM_y}
\end{flalign}
\end{subequations}
Comparison of Eqs.~(\ref{e:dCM_x}) and~(\ref{e:dCM_y}) with Eq.~(\ref{e:h4})
shows agreement through the $\mathcal{O}(\epsilon^2)$ terms.  There appears to be a discrepancy at
order $\mathcal{O}(\epsilon^3)$.  However, we have checked that this apparent discrepancy is
resolved by expanding the stochastic normal form coordinate transform to even
higher order.  For example, the seventh iteration will yield a $-4\epsilon^3
X$ term in the $y$ coordinate transform.  When added to the existing $-\epsilon^3 X$ term, there is an
agreement with the $-5\epsilon^3 x$ term in Eq.~(\ref{e:h4}).

Letting only $Y=0$ in the higher order set of equations leads to the stochastic center manifold equation.  If one
takes the expectation of this
stochastic center manifold equation and uses the following identities \cite{rob08}:
\begin{equation}
E[e^{\pm \tau}*\phi]=e^{\pm\tau}*E[\phi],
\end{equation}
\begin{equation}
E[(e^{\pm \tau}*\phi)^2]=\frac{1}{2},
\end{equation}
then one obtains the following:
\begin{flalign}
\label{e:avg_SCM}
E[y]=&X-X^3+\epsilon(-X+4X^3-3X^5)+\epsilon^2(2X-20X^3+42X^5)\nonumber\\
&+\epsilon^3(-X+16X^3-66X^5+96X^7-45X^9)\nonumber\\
&+\epsilon^4\sigma^2(-3X/2+9X^3-27X^5/2),
\end{flalign}
where the $\mathcal{O}(\epsilon^4\sigma^2)$ terms are associated with the quadratic
noise terms.

The average stochastic center manifold equation given by Eq.~(\ref{e:avg_SCM})
can now be used to solve the problem of interest.  As mentioned previously, we
are interested in computing escape rates and optimal escape paths. In
particular we have analytically found an expression for the escape rate using
the theory of large fluctuations~\cite{FW84,FM17} using both the average stochastic
center manifold and a na{\"i}ve approach wherein the noise is added after the
fact to the deterministic center manifold.
Because the noise effects occur
at such high order ($\mathcal{O}(\epsilon^4\sigma^2)$), the correction to the
na{\"i}ve approach using the 
deterministic center manifold is minimal. Details can be found in Forgoston
and Schwartz~\cite{forsch09}.

We have performed numerical computations of escape time and escape path to
compare with the analytical results. When $\sigma =0$ (no noise), the original, singularly perturbed problem given
by Eqs.~(\ref{e:Duff_sys_x})-(\ref{e:Duff_sys_y}), has three equilibrium
points given by $(-1,0)$, $(0,0)$, and
$(1,0)$.  At the initial time, $t=0$, a particle
is randomly placed near the stable,
attracting point $(1,0)$ within a circle of radius $0.1$ centered at $(1,0)$.
Equations~(\ref{e:Duff_sys_x})-(\ref{e:Duff_sys_y}) are numerically
integrated using a stochastic integrator with a constant
time step size, $\delta t$, that depends on the value of $\epsilon$,
and the time needed for the particle to escape from the basin of attraction is
determined.  This escape time is based on either the time it takes the
particle to cross the $x<-0.2$ barrier, which means the particle has escaped across the
unstable saddle, and has entered the second basin of attraction with stable,
attracting point $(-1,0)$, or when the assigned maximum time has been reached.
This computation was performed for 10\,000 particles, and the mean escape
time was determined. Figure~4 in Forgoston
and Schwartz~\cite{forsch09} demonstrates excellent agreement between the
analytical and numerical escape times.

In addition, for each of the 10\,000 particles that were initially placed in one of the
attracting basins and which later escaped from this basin, across the saddle,
and into the other basin of attraction, we retain $t=200$ worth of the
particle's path prior to escape.  By creating a histogram representing the
probability density, $p_h$, of this escape
prehistory~\cite{dmsss92}, one can see which regions of the phase space are associated with a high or low
probability of particle escape. 

\begin{figure}[h!]
\begin{center}
\includegraphics[scale=0.4]{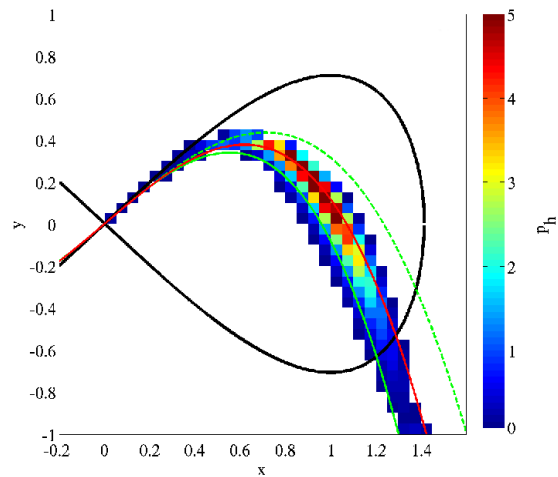}
\caption{\label{fig:hist_sig3_eps1_p10000_man}Escape path prehistory histogram
  for $\epsilon=0.1$ and $\sigma = 0.3$ overlaid with the graphs of the third-order (solid, green line), fourth-order (dashed, green
line), and fifth-order (solid, red line) center
manifold equations given by Eq.~(\ref{e:h4}). The color-bar values have been
normalized by $10^5$, and the threshold 0 value is about 9\,000. Reproduced
from Forgoston
and Schwartz~\cite{forsch09}.}
\end{center}
\end{figure}

Figure~\ref{fig:hist_sig3_eps1_p10000_man} shows a histogram of escape path
prehistory for $\epsilon=0.1$ and $\sigma = 0.3$ (so
that $D=\sigma^2 /2 = 0.045$).  The color-bar values of
Fig.~\ref{fig:hist_sig3_eps1_p10000_man} have been normalized by $10^5$.  The
threshold $0$ value in the figure is actually
about 9\,000.  Therefore, any histogram box containing less than 9\,000 events
shows up as white on the histogram. Overlaid on top of the histogram are the graphs of the
third-order, fourth-order,
and fifth-order center manifold
equations.  Each of these equations may be found by including terms of the
appropriate order from Eq.~(\ref{e:h4}). 

One can see from Fig.~\ref{fig:hist_sig3_eps1_p10000_man} that the
third-order and fourth-order center manifolds essentially bound the entire
region of escape path prehistory, while the fifth-order manifold lies along
the region of highest probability of escape.  Although it is not shown, it
should be noted that the optimal escape
path (found using large fluctuation theory) associated with the
third-order center manifold is a heteroclinic orbit from $x=\sqrt{1-\epsilon
  +2\epsilon^2}$ to $x=0$ that lies directly on top of a section of the third-order center
manifold (solid, green line in
Fig.~\ref{fig:hist_sig3_eps1_p10000_man}).  Similarly, the optimal escape
paths associated with higher order center manifolds lie directly on top of a
section of the corresponding center manifold.

Additionally, one could overlay the histogram of escape path prehistory with
the average stochastic center manifold given by Eq.~(\ref{e:avg_SCM}).
However, since the stochastic correction appears at order
$\mathcal{O}(\epsilon^4\sigma^2)$, there is no noticeable difference from the
manifolds shown in Fig.~\ref{fig:hist_sig3_eps1_p10000_man}.  Therefore, plots of the average
stochastic manifold are not shown.

Although
corrections due to noise occur at high order in the current example, other problems with quadratic nonlinearities have
noise corrections occurring at much lower order~\cite{fobisc09}, as seen in the
following example. In such cases, one must use the normal form coordinate
transform method rather than relying on the deterministic center manifold reduction.
\section{Example 2: Stochastic SEIR Epidemiological Model}\label{sec:ex2}

The interaction between deterministic and stochastic effects in population dynamics has
played, and continues to play, an important role in the modeling of
infectious diseases. The mechanistic modeling side of population dynamics
is well-known and established~\cite{Anderson91,bailey75}.  These
models typically are assumed to be useful for infinitely large, homogeneous
populations, and arise from the mean field analysis of probabilistic
models. On the other hand, when one considers finite populations, random interactions
give rise to internal noise effects, which may introduce new dynamics.  Stochastic
effects are quite prominent in finite populations, which can range from
ecological dynamics~\cite{Marion2000} to childhood epidemics in
cities~\cite{Nguyen2008,Rohani2002}.  For homogeneous populations
with seasonal forcing, noise also comes into play in the prediction
of large outbreaks~\cite{Rand1991,BillingsBS02,stolhu07}.  Specifically, external random perturbations change the probabilistic
prediction of epidemic outbreaks as well as its control~\cite{Schwartz2004}.

As a first study, we consider the Susceptible-Exposed-Infected-Recovered (SEIR)
epidemiological model with stochastic forcing.  We
could easily consider a very high-dimensional SEIR-type model where the exposed class was modeled
using hundreds of compartments.  Since the analysis is similar, we consider
the simpler standard SEIR model to demonstrate the power of the method. 

We begin by describing the stochastic version of the SEIR model found
in Schwartz and Smith~\cite{ss83}. We assume that a given population may be divided into the following four classes which evolve in time:
\begin{enumerate}
\item Susceptible class, $s(t)$, consists of those individuals who may contract the  disease.
\item Exposed class, $e(t)$, consists of those individuals who have been infected by  the disease but are not yet infectious.
\item Infectious class, $i(t)$, consists of those individuals who are capable of  transmitting the disease to susceptible individuals.
\item Recovered class, $r(t)$, consists of those individuals who are immune to the disease.
\end{enumerate}
Furthermore, we assume that the total population size, denoted as $N$, is
constant and can be normalized to $S(t)+E(t)+I(t)+R(t)=1$, where $S(t)=s(t)/N$, $E(t)=e(t)/N$,
$I(t)=i(t)/N$, and $R(t)=r(t)/N$.  Therefore, the population class variables
$S$, $E$, $I$, and $R$ represent fractions of the total population. The
governing equations for the stochastic SEIR model are 
\begin{subequations}
\begin{flalign}
&\dot{S}(t)=\mu -\beta I(t)S(t) -\mu S(t) +\sigma_1\phi_1(t),\label{e:Sdot}\\
&\dot{E}(t)=\beta I(t)S(t) -(\alpha
+\mu)E(t)+\sigma_2\phi_2(t),\label{e:Edot}\\
&\dot{I}(t)=\alpha E(t) -(\gamma +\mu)I(t)+\sigma_3\phi_3(t),\label{e:Idot}\\
&\dot{R}(t)=\gamma I(t) -\mu R(t)+\sigma_4\phi_4(t),\label{e:Rdot}
\end{flalign}
\end{subequations}
where $\sigma_i$ is the standard deviation of the noise intensity
$D_i=\sigma_i^2/2$. Each of the noise terms, $\phi_i$, describes a stochastic,
Gaussian white force that is characterized by the correlation functions given
by Eqs.~(\ref{e:mean})-(\ref{e:corr}).

Additionally, $\mu$ represents a constant birth and death rate, $\beta$ is the
contact rate, $\alpha$ is the rate of infection, so that $1/\alpha$ is the
mean latency period, and $\gamma$ is the rate of recovery, so that $1/\gamma$
is the mean infectious period. Although the contact rate $\beta$ could be
given by a time-dependent function (e.g. due to seasonal fluctuations), for
simplicity, we assume $\beta$ to be constant. Throughout this article, we use
the following parameter values: $\mu =0.02 ({\rm  year})^{-1}$, $\beta=1575.0
({\rm year})^{-1}$, $\alpha =1/0.0279 ({\rm  year})^{-1}$, and $\gamma =1/0.01
({\rm year})^{-1}$.  Disease parameters correspond to typical measles
values~\cite{ss83,bisc02}.  Note that any other biologically meaningful
parameters may be used as long as the basic reproductive rate
$R_0=\alpha\beta/[(\alpha +\mu)(\gamma +\mu)] > 1$.  The interpretation of $R_0$ is the number of secondary cases produced by a single infectious individual in
a population of susceptibles in one infectious period.

As a first approximation of stochastic effects, we have considered additive
noise.  This type of noise may result from migration into and away from the
population being considered~\cite{BvdDW08}.  Since it is difficult to estimate
fluctuating migration rates~\cite{bfg02}, it is appropriate to treat migration as an
arbitrary external noise source.  Also, fluctuations in the birth rate
manifest itself as additive noise.  Furthermore, as we are not
interested in extinction events in this article, it is not necessary to use
multiplicative noise.  In general, for the problem considered here, it is possible that a rare
event in the tail of the noise distribution may cause one or more of the $S$,
$E$, and $I$ components of the
solution to become negative.  Here, we will always assume that the
noise is sufficiently small 
so that a solution remains positive for a long
enough time to gather sufficient statistics.  Even though it is difficult to
accurately estimate the appropriate noise level from real data, our choices of
noise intensity lie within the huge confidence intervals computed in Bj{\o}rnstad
et al~\cite{bfg02}.

Although $S+E+I+R=1$ in the deterministic system, one should note that the dynamics of
the stochastic SEIR system will not necessarily have all of the components sum
to unity.  However, since the noise has zero mean, the total population will
remain close to unity on average.  Therefore, we assume that the dynamics  are
sufficiently described by Eqs.~(\ref{e:Sdot})-(\ref{e:Idot}).  It should be
noted that even if $E(t)+I(t)=0$ for some $t$, the noise allows for the
re-emergence of the epidemic.    

\subsection{Deterministic center manifold analysis}\label{sec:Dcma}
As seen in Sec.~\ref{sec:gen_theory_DCM} one can reduce the dimension of a system of equations using
deterministic center manifold theory.  
In order to make use of the center manifold theory, we transform the
deterministic version of Eqs.~(\ref{e:Sdot})-(\ref{e:Idot}) to a new system of equations that
has the necessary spectral structure.  The theory will allow us to find an invariant center manifold passing through
the fixed point to which we can restrict the transformed system.  

The transformed
evolution equations are given by 
\begin{subequations}
\begin{flalign}
\frac{dU}{d\tau}=& -\alpha_0 U+{\frac {\mu^2 \left( \gamma_0 V-\alpha_0 U \right) }{\alpha_0+\gamma_0}}-{\frac { \left( \gamma_0+\mu^2 \right) \left( \alpha_0+\mu^2 \right) \left[ \left( \alpha_0+\gamma_0 \right) U +\gamma_0 W \right] }{\alpha_0 \left( \alpha_0+\gamma_0 \right) }} \nonumber \\
& -\frac{\mu \beta}{\alpha_0+\gamma_0}  \left( \gamma_0 W+\left( \alpha_0+\gamma_0 \right) U +{\frac {\mu^2\alpha_0 \gamma_0}{ \left( \gamma_0+\mu^2 \right) \left( \alpha_0+\mu^2 \right) }} \right) \left( U+V \right) 
,\\\label{e:dU_mort}
\frac{dV}{d\tau}=& \alpha_0 U-{\frac {\mu^2 \left( \gamma_0 V -\alpha_0 U \right) }{\alpha_0+\gamma_0}}-{\frac { \left( \gamma_0+\mu^2 \right) \left( \alpha_0+\mu^2 \right) \left[ \left( \alpha_0+\gamma_0 \right) U +\gamma_0 W \right] }{\gamma_0 \left( \alpha_0+\gamma_0 \right) }} \nonumber \\
& -\frac{\mu \beta \alpha_0}{ \gamma_0 \left( \alpha_0+\gamma_0 \right)} \left( \gamma_0 W+\left( \alpha_0+\gamma_0 \right) U +{\frac {\mu^2\alpha_0 \gamma_0}{ \left( \gamma_0+\mu^2 \right) \left( \alpha_0+\mu^2 \right) }} \right) \left( U+V \right),\\\label{e:dV_mort}
\frac{dW}{d\tau}=& -\alpha_0 U - \left( \gamma_0+\mu^2 \right) \left( U+W \right) +{\frac { \left( \gamma_0+\mu^2 \right) \left( \alpha_0+\mu^2 \right) \left[ \left( \alpha_0+\gamma_0 \right) U +\gamma_0 W \right] }{\alpha_0 \gamma_0}} \nonumber \\
& -\mu^2V + \frac{\mu \beta}{\gamma_0} \left( \gamma_0 W+\left( \alpha_0+\gamma_0 \right) U +{\frac {\mu^2\alpha_0 \gamma_0}{ \left( \gamma_0+\mu^2 \right) \left( \alpha_0+\mu^2 \right) }} \right) \left( U+V \right)
 ,\\\label{e:dW_mort}
\frac{d\mu}{d\tau}=&0,\\\label{e:dmu_mort} 
\end{flalign}
\end{subequations}
where 
\begin{subequations}
\begin{flalign}
&U = \frac{-\gamma_0}{\alpha_0 + \gamma_0}\bar{E},\label{e:U}\\
&V = \bar{S} + \frac{\gamma_0}{\alpha_0 + \gamma_0}\bar{E},\label{e:V}\\
&W = \bar{I} + \bar{E},\label{e:W}
\end{flalign}
\end{subequations}
$\bar{S}(t)=S(t)-S_0$, $\bar{E}(t)=E(t)-E_0$, and $\bar{I}(t)=I(t)-I_0$, 
\begin{equation}
\label{e:fixed_pt_mort}
(S_0,E_0,I_0)=\left (
{\frac {  \left( \gamma+\mu   \right)  \left( \alpha+\mu  \right) }{\beta \alpha}},
{\frac {\mu   }{\alpha+\mu }}-{\frac {\mu   \left( \gamma+\mu  \right) }{\alpha \beta}} ,
{\frac {\mu \alpha}{ \left( \gamma+\mu  \right)  \left( \alpha+\mu  \right) }}-{\frac {\mu  }{\beta}}
  \right )
\end{equation}
is a fixed point corresponding to the endemic state, and $\alpha=\alpha_0/\mu$ and $\gamma=\gamma_0/\mu$,
where $\alpha_0$ and $\gamma_0$ are $\mathcal{O}(1)$.
Details
regarding the transformation can be found in Forgoston et al~\cite{fobisc09}.

The Jacobian of Eqs.~(\ref{e:dU_mort})-(\ref{e:dmu_mort}) to zeroth-order in
$\mu$ and evaluated at the origin is
\begin{equation}
\left [ \begin{array}{cc|ccccccccc}
-(\alpha_0+\gamma_0) & & & &0 & & & -\frac{\gamma_0^2}{(\alpha_0+\gamma_0)}&
& & 0 \\
& & & & & & & & & &\\
\hline
& & & & & & & & & &\\
0 & & & &0 & & & -\frac{\alpha_0\gamma_0}{(\alpha_0+\gamma_0)}& & & 0 \\
& & & & & & & & & &\\
0 & & & &0 & & & 0 & & & 0 \\
& & & & & & & & & &\\
0 & & & &0 & & & 0 & & & 0 \end{array} \right ],
\end{equation} 
which shows that Eqs.~(\ref{e:dU_mort})-(\ref{e:dmu_mort}) may be rewritten in the form
\begin{flalign}
&\frac{d{\bf x}}{d\tau}= {\bf A}{\bf x} + {\bf f}({\bf x},{\bf  y},\mu),\label{e:CMformx}\\
&\frac{d{\bf y}}{d\tau}= {\bf B}{\bf y} + {\bf g}({\bf x},{\bf  y},\mu),\label{e:CMformy}\\
&\frac{d\mu}{d\tau}=0,
\end{flalign}
where ${\bf x}=(U)$, ${\bf y}=(V,W)$, ${\bf A}$ is a constant matrix with
eigenvalues that have negative real parts, ${\bf B}$ is a constant matrix with
eigenvalues that have zero real parts, and ${\bf f}$ and ${\bf g}$ are
nonlinear functions in ${\bf x}$, ${\bf y}$ and $\mu$.  In particular,
\begin{equation}
{\bf A}=\left [ \begin{array}{c}
-(\alpha_0 +\gamma_0)\end{array} \right ], \,\,\,\,\,\,\,\,\,\,\,\,\,\,\, {\bf B}=\left [
\begin{array}{cccc}
0 &&& -\frac{\alpha_0\gamma_0}{(\alpha_0 +\gamma_0)}\\
&&&\\
0 &&& 0\end{array} \right ].
\end{equation}

 Therefore, the system will rapidly collapse onto a lower-dimensional manifold
 given by center manifold theory~\cite{car81}. Furthermore, we know that the center manifold is given by
\begin{equation}
\label{e:CM}
U=h(V,W,\mu),
\end{equation}
where $h$ is an unknown function.

Substitution of Eq.~(\ref{e:CM}) into Eq.~(\ref{e:dU_mort}) leads to the
following center manifold condition:
\begin{flalign}
\label{e:CMcond}
&\frac{\partial h(V,W,\mu)}{\partial V}\frac{dV}{d\tau}+\frac{\partial
  h(V,W,\mu)}{\partial W}\frac{dW}{d\tau}= -\alpha_0 h(V,W,\mu)+ {\frac {\mu^2
    \left[ \gamma_0 V-\alpha_0 h(V,W,\mu) \right]
  }{\alpha_0+\gamma_0}}-&&\nonumber\\
&{\frac { \left( \gamma_0+\mu^2 \right) \left( \alpha_0+\mu^2 \right) \left[ \left( \alpha_0+\gamma_0 \right) h(V,W,\mu) +\gamma_0 W \right] }{\alpha_0 \left( \alpha_0+\gamma_0 \right) }}- \nonumber \\
& \frac{\mu \beta}{\alpha_0+\gamma_0}  \left( \gamma_0 W+\left( \alpha_0+\gamma_0 \right) h(V,W,\mu) +{\frac {\mu^2\alpha_0 \gamma_0}{ \left( \gamma_0+\mu^2 \right) \left( \alpha_0+\mu^2 \right) }} \right) \left( h(V,W,\mu)+V \right).
\end{flalign}
In general, it is not possible to solve the center manifold condition for the unknown function, $h(V,W,\mu)$. Therefore, we perform the following Taylor series expansion of $h(V,W,\mu)$ in $V$, $W$, and $\mu$:
\begin{flalign}
\label{e:TS}
h(V,W,\mu)=&h_0 + h_2 V +h_3 W +h_\mu\mu + h_{22}V^2 + h_{23}VW +h_{33}W^2
+\nonumber\\
&h_{\mu 2}\mu V + h_{\mu 3}\mu W + h_{\mu\mu}\mu^2 + \ldots,
\end{flalign} 
where $h_0$, $h_2$, $h_3$, $h_\mu$, $\ldots$ are unknown coefficients that are found by substituting the Taylor series expansion into the center manifold condition and equating terms of the same order. By carrying out this procedure using a second-order Taylor series expansion of $h$, the center manifold equation is
\begin{equation}
\label{e:CMeq}
U=-\frac{\gamma_0^2}{(\alpha_0+\gamma_0)^2}W + \mathcal{O}(\epsilon^3),
\end{equation}
where $\epsilon =|(V,W,\mu)|$ so that $\epsilon$ provides a count of the
number of $V$, $W$, and $\mu$ factors in any one term.  Substitution of Eq.~(\ref{e:CMeq}) into Eqs.~(\ref{e:dV_mort}) and~(\ref{e:dW_mort}) leads to the following reduced system of evolution equations which describe the dynamics on the center manifold: 
\begin{subequations}
\begin{flalign}
\frac{dV}{d\tau}=& -{\frac {\mu^2{\gamma_0}^2\alpha_0  W}{ \left(
      \alpha_0+\gamma_0 \right) ^{3}}} -{\frac {  \mu^{4}\alpha_0 W}{ \left(
      \alpha_0+\gamma_0 \right)^2}}  -{\frac {\gamma_0 \mu^2
    V}{\alpha_0+\gamma_0}}-\nonumber\\
&{\frac { \left( \gamma_0+\mu^2 \right)  \alpha_0 W}{\alpha_0+\gamma_0}}-
\frac{\beta {\alpha_0}^2\mu}{ \left( \alpha_0+\gamma_0
  \right)^2} \left(  W+{\frac {\mu^2 \left( \alpha_0+\gamma_0 \right) }{ \left( \gamma_0+\mu^2 \right)  \left( \alpha_0+\mu^2 \right) }} \right)  \left(  V-{\frac {{\gamma_0}^2 W}{ \left( \alpha_0+\gamma_0 \right)^2}} \right) 
,\label{e:dV_mort_red}\\
\frac{dW}{d\tau}=& {\frac {\mu^2 {\gamma_0}^2W }{ \left( \alpha_0+\gamma_0
    \right)^2}}+{\frac {\mu^{4} W}{\alpha_0+\gamma_0}}-\mu^2 V+\nonumber\\
& \frac{\beta
  \mu \alpha_0}{ \alpha_0+\gamma_0 }  \left(  W+{\frac {\mu^2 \left( \alpha_0+\gamma_0 \right) }{ \left( \gamma_0+\mu^2 \right)  \left( \alpha_0+\mu^2 \right) }} \right)  \left(  V-{\frac {{\gamma_0}^2 W}{ \left( \alpha_0+\gamma_0 \right)^2}} \right) 
.\label{e:dW_mort_red} 
\end{flalign}
\end{subequations}

\subsection{Projection of the noise onto the stochastic center
 manifold}\label{sec:Incorrect_proj}

The stochastic SEIR system given by Eqs.~(\ref{e:Sdot})-(\ref{e:Idot}) may be
transformed in a manner similar to what was done in the preceding section for
the deterministic form of the system. The transformation will have an effect
on the noise terms. The transformed stochastic terms are still additive,
Gaussian noise processes.  However, the transformation has
acted on the original stochastic terms $\phi_1$, $\phi_2$, and $\phi_3$ to
create new noise processes which have a variance different from that of the
original noise processes. This transformed system of equations is an exact transformation of the system of equations given by
Eqs.~(\ref{e:Sdot})-(\ref{e:Idot}) and can be found in Forgoston et al~\cite{fobisc09}.

\subsubsection{Reduction of the stochastic SEIR model}\label{sec:ipcm_mort}

It is tempting to consider the reduced stochastic model found by substitution
of Eq.~(\ref{e:CMeq}) into the tranformed equations, so that one has the following stochastic evolution equations (that hopefully describe the dynamics on the stochastic center manifold):
\begin{subequations}
\begin{flalign}
\frac{dV}{d\tau}=& -{\frac {\mu^2{\gamma_0}^2\alpha_0  W}{ \left( \alpha_0+\gamma_0 \right) ^{3}}} -{\frac {  \mu^{4}\alpha_0 W}{ \left( \alpha_0+\gamma_0 \right)^2}}  -{\frac {\gamma_0 \mu^2 V}{\alpha_0+\gamma_0}}-{\frac { \left( \gamma_0+\mu^2 \right)  W \alpha_0}{\alpha_0+\gamma_0}}- \nonumber \\
& \frac{\beta {\alpha_0}^2\mu}{ \left( \alpha_0+\gamma_0 \right)^2}  \left(  W+{\frac {\mu^2 \left( \alpha_0+\gamma_0 \right) }{ \left( \gamma_0+\mu^2 \right)  \left( \alpha_0+\mu^2 \right) }} \right)  \left(  V-{\frac {{\gamma_0}^2 W}{ \left( \alpha_0+\gamma_0 \right)^2}} \right) 
+\sigma_5\phi_5,\label{e:dV_mort_red_stoch}\\
\frac{dW}{d\tau}=& {\frac {\mu^2 {\gamma_0}^2W }{ \left( \alpha_0+\gamma_0 \right)^2}}+{\frac {\mu^{4} W}{\alpha_0+\gamma_0}}-\mu^2 V+  \nonumber \\
&  \frac{\beta \mu \alpha_0}{ \alpha_0+\gamma_0 }  \left(  W+{\frac {\mu^2 \left( \alpha_0+\gamma_0 \right) }{ \left( \gamma_0+\mu^2 \right)  \left( \alpha_0+\mu^2 \right) }} \right)  \left(  V-{\frac {{\gamma_0}^2 W}{ \left( \alpha_0+\gamma_0 \right)^2}} \right) + \sigma_6\phi_6,\label{e:dW_mort_red_stoch} 
\end{flalign}
\end{subequations}
where
\begin{subequations}
\begin{flalign}
\sigma_5\phi_5 = & \mu\sigma_1\phi_1 + \frac{\mu\gamma_0}{\alpha_0
  +\gamma_0}\sigma_2\phi_2,\label{e:sig5phi5}\\
\sigma_6\phi_6 = & \mu\sigma_2\phi_2 +\mu\sigma_3\phi_3\label{e:sig6phi6}
\end{flalign}
\end{subequations}
are transformed stochastic terms.

One should note that Eqs.~(\ref{e:dV_mort_red_stoch})-(\ref{e:dW_mort_red_stoch}) also can be found by na\"{i}vely adding the
stochastic terms to the reduced system of evolution equations for the
deterministic problem (Eqs.~(\ref{e:dV_mort_red})-(\ref{e:dW_mort_red})).  This type of stochastic center manifold reduction has been done for the case
of discrete noise~\cite{bisc02}.  Additionally, in many other fields
(e.g. oceanography, solid mechanics, fluid mechanics), researchers have
performed stochastic model reduction using a Karhunen-Lo\`{e}ve expansion
(principal component analysis, proper orthogonal decomposition)~\cite{doghre07,vewaka08}.
However, this linear projection does not properly capture the nonlinear
effects.  Furthermore, one must subjectively choose the number of modes needed
for the expansion.  Therefore, even though the solution to the reduced model
found using this technique may have the correct statistics, individual
solution realisations will not agree with the original, complete solution.  

\begin{figure}[t!]
\centering
\includegraphics[width=8.5cm,height=5cm]{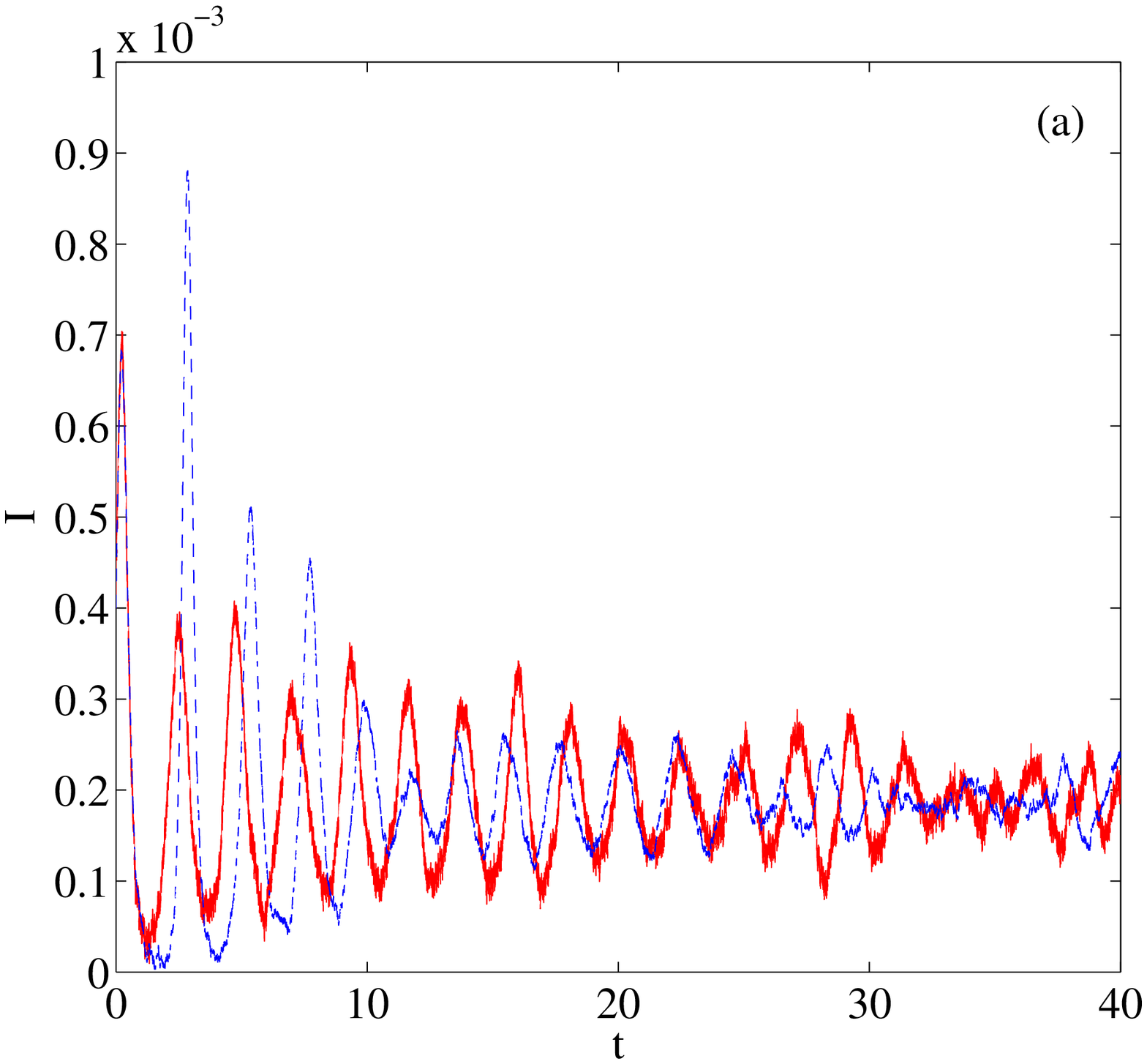}\\
\includegraphics[width=8.5cm,height=5cm]{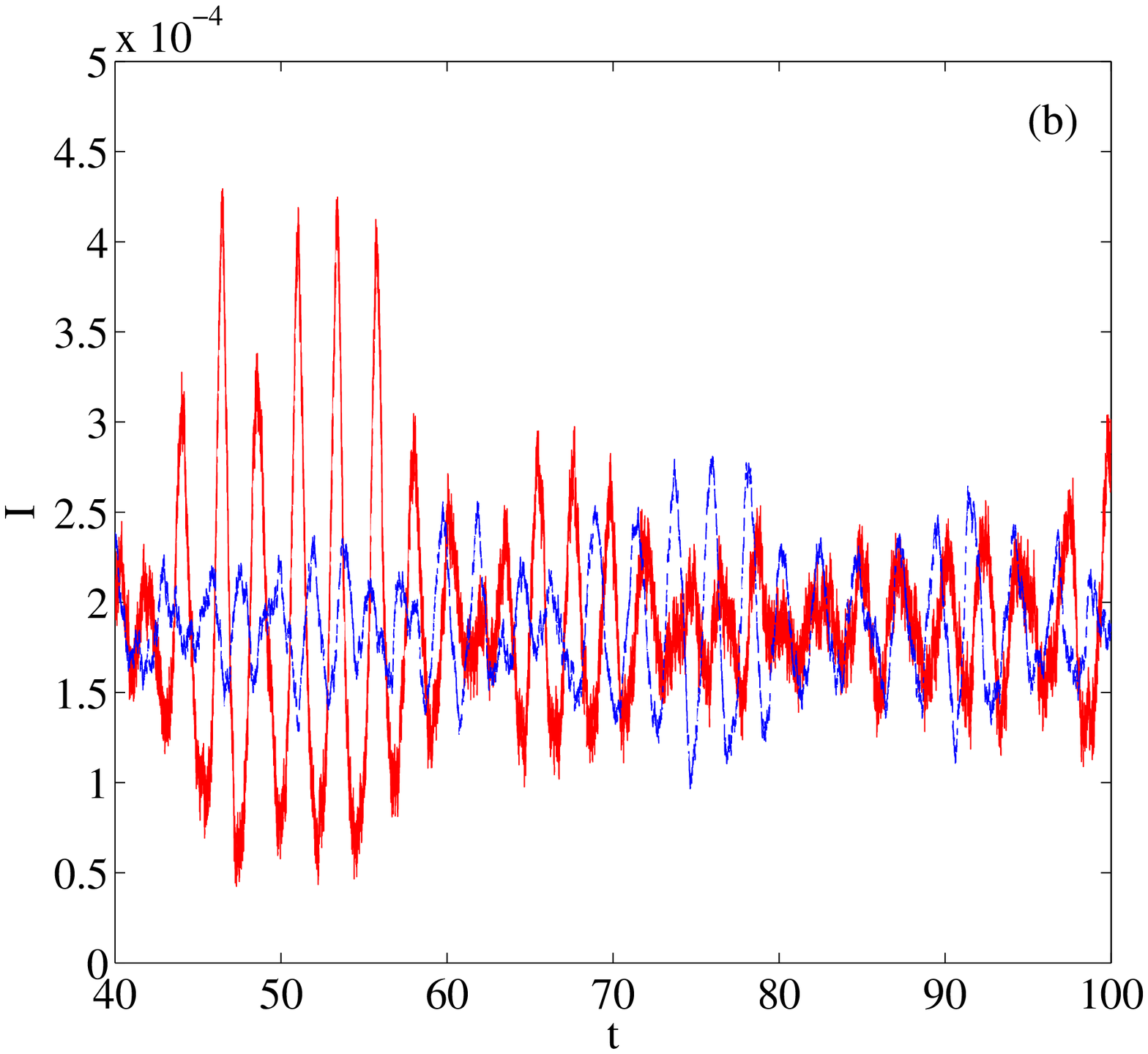}
\caption{\label{fig:TS_full_reducedCM} Time series of the fraction of the
  population that is  infected with a disease, $I$, computed using the
  complete, stochastic system of transformed equations of the SEIR model (red, solid line),
  and computed using the reduced system of equations of the SEIR model that is based on the
  deterministic center manifold with a replacement of the noise terms
  (Eqs.~(\ref{e:dV_mort_red_stoch})-(\ref{e:dW_mort_red_stoch})) (blue,
  dashed line). The standard deviation of the noise intensity used in the
  simulation is $\sigma_i =0.0005$, $i=4,5,6$.  The time series is shown for (a)
  $t=0$ to $t=40$, and for (b) $t=40$ to $t=100$. Reproduced
from Forgoston et al~\cite{fobisc09}.}
\end{figure}

It can be seen in Fig.~\ref{fig:TS_full_reducedCM} that Eqs.~(\ref{e:dV_mort_red_stoch})-(\ref{e:dW_mort_red_stoch}) do not contain the correct projection of the
noise onto the center manifold.  Therefore, when solving the reduced system,  one does not obtain the correct solution. Such errors in stochastic  epidemic modeling impact the prediction of disease outbreak when modeling  the spread of a disease in a population. 

Figures~\ref{fig:TS_full_reducedCM}(a)-(b) compares the time series of the fraction of
the population that is infected with a disease, $I$, computed using the
complete, stochastic system of transformed equations of the SEIR model
 with the time series of
$I$ computed using the reduced system of equations of the SEIR model that is based on the
deterministic center manifold with a replacement of the noise terms
(Eqs.~(\ref{e:dV_mort_red_stoch}) and~(\ref{e:dW_mort_red_stoch})).
Figure~\ref{fig:TS_full_reducedCM}(a) shows the initial transients, while
Fig.~\ref{fig:TS_full_reducedCM}(b) shows the time series after the transients have
decayed.  One can see that the solution computed using the reduced system quickly
becomes out of phase with the solution of the complete system.  Use of this
reduced system would lead to an incorrect prediction of the initial disease
outbreak.  Additionally, the predicted amplitude of the initial outbreak would
be incorrect.  The poor agreement, both in phase and amplitude, between the
two solutions continues for long periods of time as seen in Fig.~\ref{fig:TS_full_reducedCM}(b).  We also have
computed the cross-correlation of the two time series shown in
Fig.~\ref{fig:TS_full_reducedCM}(a)-(b) to be approximately 0.34.  Since the cross-correlation measures the
similarity between the two time series, this low value quantitatively suggests
poor agreement between the two solutions.

\subsubsection{Correct projection of the noise onto the stochastic center
 manifold}\label{sec:correct_proj}
To project the noise correctly onto the center manifold, we derive a
normal form coordinate transform using the principles of Sec.~\ref{sec:gen_theory_SCM} for the complete, stochastic system of
transformed equations of the SEIR model. 

The stochastic system of equations (which satisfies the
necessary spectral requirements) in $\left (U,V,W\right )^T$ coordinates is
transformed to a new $\left (Y,X_1,X_2\right )^T$ coordinate system using a
near-identity stochastic coordinate transform given as
\begin{subequations}
\begin{flalign}
U=&Y+ \xi\left
 (Y,X_1,X_2,\tau\right ),\\
V=&X_1+ \eta\left (Y,X_1,X_2,\tau\right ),\\
W=&X_2+ \rho\left (Y,X_1,X_2,\tau\right ),
\end{flalign}
\end{subequations}
where the specific form of $\xi\left (Y,X_1,X_2,\tau\right )$, $\eta\left
  (Y,X_1,X_2,\tau\right )$, and $\rho\left (Y,X_1,X_2,\tau\right )$ is found using an iterative procedure.  

Several iterations lead to coordinate transforms for $U$, $V$, and $W$ along
with evolution equations describing the $Y$-dynamics, $X_1$-dynamics, and
$X_2$-dynamics in the new coordinate system. The $Y$-dynamics have exponential
decay to the $Y=0$ slow manifold. Substitution of $Y=0$ leads to complicated
expressions for the coordinate transforms - they can be found in Forgoston et al~\cite{fobisc09}.

We are interested in the long-term slow processes. Since the memory integrals
are fast-time processes, we neglect them along with the higher-order
multiplicative terms to obtain the following expressions for $U$, $V$, and $W$:
\begin{subequations}
\begin{flalign}
U = &  - \frac{\gamma_0^2 X_2}{\left (\alpha_0 + \gamma_0\right )^2}
-  \frac{\mu\beta X_1}{\left (\alpha_0 + \gamma_0\right )} \left(
\frac{\mu^2}{\left (\alpha_0 + \gamma_0\right )}
+\frac{\gamma_0 X_2}{\left (\alpha_0 + \gamma_0\right )}
+ \mu^2 X_2 \right)&& \nonumber\\
&+\frac{\mu^2\gamma_0}{\left (\alpha_0+\gamma_0\right )}\left( X_1-\frac{2X_2}{\alpha_0} \right)
,\label{e:U_cm_mort}\\
V = & ~ X_1,&&\\
W = & ~ X_2.&&
\end{flalign}
\end{subequations}
Note that Eq.~(\ref{e:U_cm_mort}) is the deterministic center manifold
equation, and at first-order, matches the center manifold equation that was
found previously (Eq.~(\ref{e:CMeq})).

Substitution of $Y=0$ and neglecting all multiplicative noise terms and memory
integrals using the argument from above (so that we consider only first-order noise terms) leads to the
following reduced system of evolution equations on the center manifold:
\begin{subequations}
\begin{flalign}
\frac{dX_1}{d\tau} = F(X_1(\tau),X_2(\tau)),\label{e:F}\\
\frac{dX_2}{d\tau} = G(X_1(\tau),X_2(\tau))\label{e:G}.
\end{flalign}
\end{subequations}
The specific form of $F$ and $G$ in Eqs.~(\ref{e:F}) and~(\ref{e:G}) are
complicated, and can be found in Forgoston et al~\cite{fobisc09}.

\begin{figure}[t!]
\centering
\includegraphics[width=8.5cm,height=5cm]{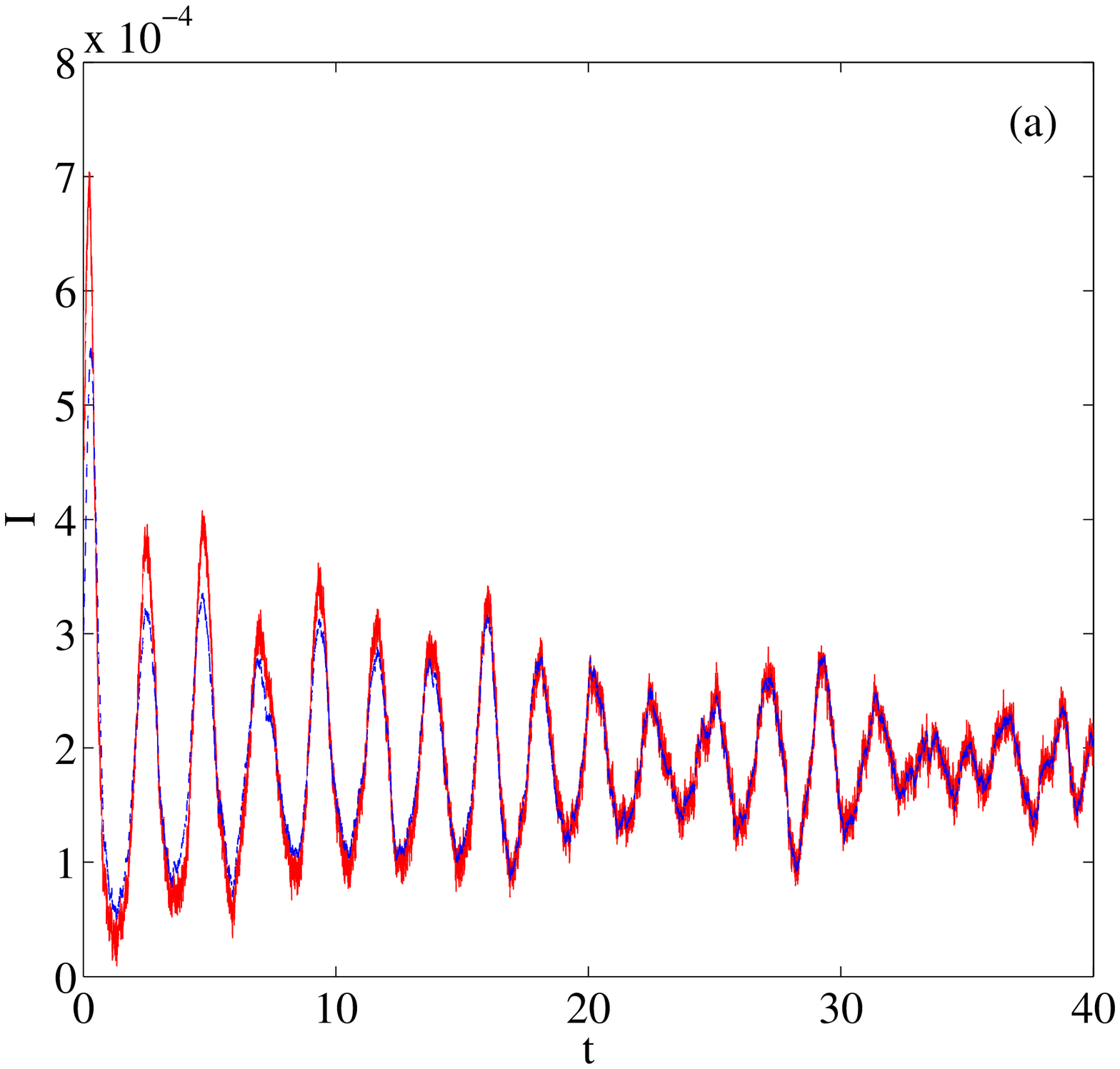}\\
\includegraphics[width=8.5cm,height=5cm]{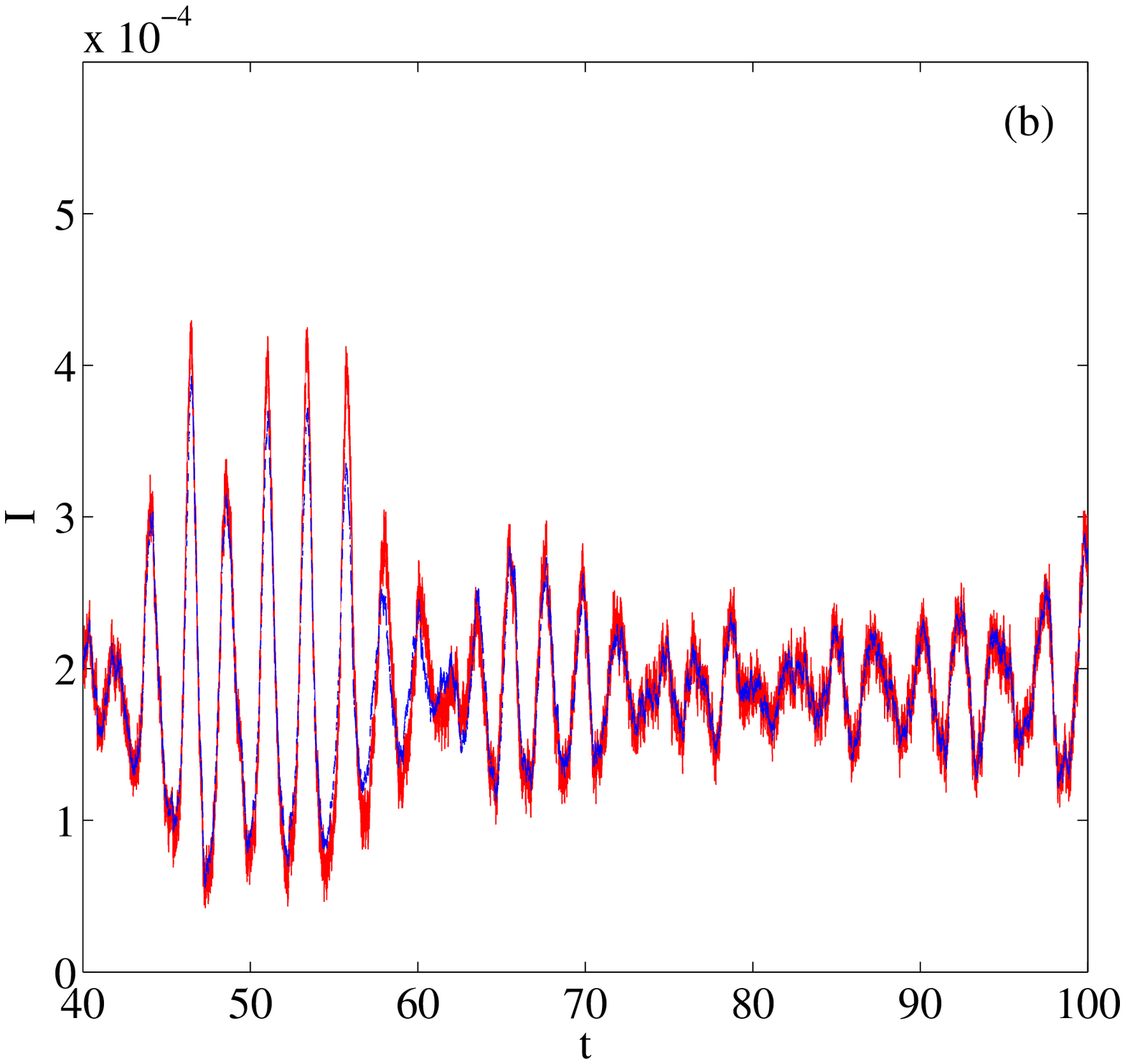}
\caption{\label{fig:TS_full_reducedNF} Time series of the
  fraction of the population that is  infected with a disease, $I$, computed
  using the complete, stochastic system of transformed equations of the SEIR
  model (red, solid
  line), and computed using the reduced system of equations of the SEIR model
  that is found using the stochastic normal form coordinate transform
  (Eqs.~(\ref{e:F})-(\ref{e:G})) (blue, dashed line). The standard deviation of
  the noise intensity used in the simulation is $\sigma_i =0.0005$,
  $i=4,5,6$.  The time series is shown for (a)
  $t=0$ to $t=40$, and for (b) $t=40$ to $t=100$. Reproduced
from Forgoston et al~\cite{fobisc09}.}
\end{figure}

Figures~\ref{fig:TS_full_reducedNF}(a)-(b) compares the time series of the fraction of
the population that is infected with a disease, $I$, computed using the
complete, stochastic system of transformed equations of the SEIR model
 with the time series of
$I$ computed using the reduced system of equations of the SEIR model that is found using
the stochastic normal form coordinate transform (Eqs.~(\ref{e:F})-(\ref{e:G})).
Figure~\ref{fig:TS_full_reducedNF}(a) shows the initial transients, while
Fig.~\ref{fig:TS_full_reducedNF}(b) shows the time series after the transients
have decayed.  One can see that there is excellent agreement between the two
solutions.  The initial outbreak is successfully captured by the reduced
system.  Furthermore, Fig.~\ref{fig:TS_full_reducedNF}(b) shows that the reduced system accurately predicts recurrent
outbreaks for a time scale that is orders of magnitude longer than the
relaxation time.  This is not surprising since the solution decays
exponentially throughout the transient and then remains close to the
lower-dimensional center manifold.  Since we are not looking at periodic
orbits, there are no secular terms in the asymptotic expansion, and the result
is valid for all time.  Additionally, any noise drift on the center manifold
results in bounded solutions due to sufficient dissipation transverse to the manifold.  The cross-correlation of the two time series shown in
Fig.~\ref{fig:TS_full_reducedNF} is approximately 0.98, which quantitatively suggests
there is excellent agreement between the two solutions.

Unlike example 1, where the average stochastic center manifold and deterministic center
manifold are virtually identical,  the SEIR model considered in example 2 has terms at low order in the normal form transform which cause
a significant difference between the average stochastic center manifold and
the deterministic manifold.  Therefore, as we have demonstrated, when working
with the SEIR model, one must use
the normal form coordinate transform method to correctly project the noise
onto the center manifold.

\section{Conclusions}\label{sec:conc}

This review highlights some of our results in stochastic model reduction that use the
normal form coordinate transform method of A.J. Roberts. It is increasingly
apparent that small external effects can play a significant role in the
dynamics of models used for real world applications. Therefore, researchers
are revisiting existing methods and developing new ones to quantify
noise-induced phenomena in stochastic models. Time series prediction is a
classic problem in deterministic dynamical systems. As we develop methods to
approach its counterpart in stochastic models, we recognize the limitations of
high dimensional analysis. Moreover, many techniques that are
  probabilistic in nature for stochastic systems can be carried over to
  deterministic complex systems, such as chaos or turbulence.  Reduced models provide tractable systems to analyze
and test, both analytically and computationally. However, there is still much
work to be done in moving to high-dimensional systems, including dynamics on
large-scale heterogeneous networks where mean field analysis typically
fails~\cite{schwartz2016epidemic}, and truly infinite stochastic systems arising
from delay~\cite{schwartz2015noise}, where the force of the noise needs to be
known in the past and future.

In this work, we demonstrate
the importance of properly projecting the noise source through the
reduction. The method of Roberts is an important contribution as it not only
enables one to reduce the dimensionality of the dynamics, but also provides a
way to separate all slow processes from all fast processes. We have
demonstrated in these two examples how deriving the deterministic manifold may
not be sufficient, and in this case, how one can employ the normal form
coordinate transform method to properly capture the stochastic manifold. The result in the SEIR example is a closed form system that can be analyzed and tested for extended time series prediction. Stochastic model reduction is a general approach that can be used to analyze the dynamics of generic stochastic models with separated time scales.

\section*{Acknowledgments}
EF and LB are supported by the National Science Foundation award DMS
\#1418956. IBS is funded under the Office of Naval Research award N0001412WX20083, and the NRL Base Research Program N0001412WX30002.

We thank Tony (A.J. Roberts) for his conversations, comments and insight over the years
regarding the work presented in this article.
\bibliographystyle{ws-rv-van}
\bibliography{Refs_final}

\end{document}